\documentclass[11pt,reqno]{amsart}
\usepackage{hyperref}
\usepackage{url}
\usepackage{amssymb}
\usepackage[T1]{fontenc}
\usepackage{lmodern}
\usepackage[utf8]{inputenc}
\usepackage{ifthen}
\usepackage{mathtools}

\newtheorem{theorem}{Theorem}

\theoremstyle{definition}

\theoremstyle{remark}

\numberwithin{equation}{section}

\usepackage{tikz,enumerate}
\allowdisplaybreaks[4]
\begin{document}

\title{On the difference of partial theta functions}

\author{Alexander Berkovich}

%
%

%


%
\begin{abstract}
Sums of the form $\sum_{n\geq 0} (-1)^nq^{(n-1)n/2} x^n $ are called partial theta functions. In his lost
notebook, Ramanujan recorded many identities for those functions. In 2003, Warnaar found an elegant formula for a sum of two partial theta functions. Subsequently, Andrews and Warnaar established a similar result for the product of two partial theta functions. In this note, I discuss the relation between the Andrews-Warnaar identity and the (1986) product formula due to Gasper and Rahman. I employ nonterminating extension of Sears-Carlitz transformation for $\, _3\phi_2$ to provide a new elegant proof for a companion identity for the difference of two partial theta series. This difference formula first appeared in the work of Schilling-Warnaar (2002). Finally, I show that Schilling-Warnnar (2002) and Warnaar (2003) formulas are, in fact, equivalent.

\end{abstract}

\maketitle

\section{Introduction}

I begin by recalling some standard notations and terminology in \cite{Gasper-Rahman}. Let $a, q$ be complex numbers with $0<|q|<1$. Then the $q$-shifted factorial is defined by
$$(a)_n = (a;q)_n :=\begin{cases} \prod_{j=0}^{n-1}(1-aq^j)\,\,\text{if $n>0$,}\\
1 \,\,\text{if $n=0,$}\end{cases}$$
for $n\in \mathbf N$. For the sake of brevity, I often use
$$(a_1,a_2,\cdots, a_m;q)_n =(a_1,a_2,\cdots , a_m)_n := (a_1)_n(a_2)_n\cdots (a_m)_n.$$
The basic hypergeometric series is defined by
$$\, _{r+1}\phi_r\left(\begin{matrix*}[l] a_1,a_2,\cdots, a_{r+1};q,z \\
b_1,b_2,\cdots , b_r\end{matrix*}\right) := \sum_{j=0}^\infty \dfrac{(a_1,a_2,\cdots , a_{r+1})_j}{(q,b_1,b_2,\cdots , b_r)_j} z^j.$$
The celebrated Jacobi's triple product identity can be written as
\begin{equation}\label{eq1}
\sum_{n=-\infty}^\infty (-1)^n q^{n(n-1)/2} x^n = (q,x,q/x)_\infty,
\end{equation}
and
$$(a_1,a_2,\cdots, a_m)_\infty :=\lim_{n\rightarrow \infty} (a_1,a_2,\cdots , a_m)_n.$$
Another famous result of Jacobi
 \begin{equation}\label{eq2}
\sum_{n=1}^\infty (-1)^n q^{n(n-1)/2}(2n-1)  = -(q)_\infty^3
\end{equation}
is an immediate corollary of \eqref{eq1}.

In \cite{Warnaar}, Warnaar proved the following generalization of \eqref{eq1}
\begin{equation}\label{eq3}
1+\sum_{n\geq 1} (-1)^n q^{n(n-1)/2} (x^n+y^n) =(q,x,y)_\infty \sum_{n\geq 0}\dfrac{(xy/q)_{2n} q^n}{(q,x,y,xy)_n}.
\end{equation}
Indeed, it is easy to check that \eqref{eq1} follows from \eqref{eq3} upon setting $y=q/x$ and noting that
$$\sum_{n\geq 0}\dfrac{(1)_{2n}q^n}{(q,x,q/x,q)_n}=1,$$
thanks to $(1)_n=\delta_{n,0}.$

Sums of the form $$\sum_{n\geq 1} (-1)^n q^{n(n-1)/2} x^n$$ are called partial theta series owing to the fact that
$$\sum_{n=-\infty}^\infty (-1)^nq^{n(n-1)/2} x^n$$ is referred to as a complete theta function or just
a theta function. Partial theta functions play a prominent role in Ramanujan's Lost Notebook \cite{LNB}. True
to himself, Ramanujan gave no proof of any of the partial theta function formulae in \cite{LNB}, making
it virtually impossible to determine how he discovered them. Ramanujan's identities were explicated by Andrews
in \cite{Andrews}. For more recent work motivated by \cite{LNB} the reader is invited to examine \cite{Kang} and
references there. Before I commence my calculations let me pause for a light but thoughtful
\medskip

\noindent {\it Remark.} The term ``Lost Notebook'' is somewhat imprecise, because Ramanujan's work has never been lost. However, this romantic misnome inspired a generation of young mathematicians throughout the world. It made
them feel like Indiana Jones on the quest for the Lost Ark of the Convenant. And so, let me retain this
spacious term. After all, an asteroid is not a star-like objet as its name may suggest.

\medskip

The remainder of this note is organized as follows. In the next section, I discuss sum, product and difference
formulae for partial theta series. In particular, I prove that
\begin{equation}\label{eq4} \sum_{n\geq 1} (-1)^n q^{n(n-1)/2} \dfrac{x^n-y^n}{x-y}=-(q,qx,qy)_\infty \sum_{n\geq 0}\dfrac{(xy)_{2n} q^n}{(q,qx,qy,xy)_n}.\end{equation}
I remark that this formula is equivalent to the Lemma 4.3 in \cite{SW}. Also, in the next section I prove the equivalence of \eqref{eq3} and \eqref{eq4}. Section 3 contains some appealing corollaries of \eqref{eq4}. I conclude with an interesting octonic transformation formula.

\section{The middle section}

I start by showing that Heine's first transformation formula \cite[(III.1)]{Gasper-Rahman}
\begin{equation}\label{eq5}
\, _2\phi_1\left(\begin{matrix*}[l] a,b \\ c \end{matrix*};q,z\right)
=\dfrac{(b,az)_\infty}{(c,z)_\infty} \, _2\phi_1\left(\begin{matrix*}[l] c/b, z \\ az  \end{matrix*};q,b\right)\end{equation}
can be employed to find another useful representation for the partial theta function. Indeed,
\begin{align}
\notag
\sum_{n\geq 0}(-1)^n q^{n(n-1)/2} x^n
 =\lim_{\rho \rightarrow \infty} \, _2\phi_1\left(\begin{matrix*}[l] \rho, q \\ 1/\rho \end{matrix*};q,x/\rho\right)&=
 (q,x)_\infty \, _2\phi_1\left(\begin{matrix*}[l] 0, 0 \\ x \end{matrix*};q,q\right) \\
 &=(q,x)_\infty \sum_{n\geq 0}\dfrac{q^n}{(q,x)_n}.\label{eq6}\end{align}

In \cite{Andrews-Warnaar}, Andrews and Warnaar proved that
\begin{equation}\label{eq7}
\sum_{n\geq 0}(-1)^n q^{n(n-1)/2} x^n \sum_{m\geq 0} (-1)^m q^{m(m-1)/2} y^m
=(q,x,y)_\infty \sum_{n\geq 0}\dfrac{(xy/q)_{2n} q^n}{(q,x,y,xy/q)_n}.
\end{equation}
Moreover, they showed that \eqref{eq3} is an immediate corollary of \eqref{eq7}.

It turns out that \eqref{eq7} is a special case of the Gasper-Rahman product formula
\cite[(8.8.18)]{Gasper-Rahman}
\begin{align}
\notag
&\, _2\phi_1\left(\begin{matrix*}[l] a, b \\ c \end{matrix*};q,z\right)
\, _2\phi_1\left(\begin{matrix*}[l] a, aq/c \\ aq/b \end{matrix*};q,z\right)
\\ \notag
&=\dfrac{(az,abz/c)_\infty}{(z,bz/c)_\infty} \, _6\phi_5\left(\begin{matrix*}[l] a, c/b, \sqrt{ac/b},
 - \sqrt{ac/b}, \sqrt{acq/b},- \sqrt{acq/b}\\  aq/b, c, ac/b, az, cq/bz \end{matrix*};q,q\right)\\
 &+\dfrac{(a,c/b,az,bz,azq/c)_\infty}{(c,aq/b,z,z,c/bz)_\infty}
 \, _6\phi_5\left(\begin{matrix*}[l] z, abz/c, z\sqrt{ab/c},
 -z\sqrt{ab/c}, z\sqrt{abq/c},-z\sqrt{abq/c}\\  az, bz, azq/c, bzq/c, abz^2/c \end{matrix*};q,q\right)\label{eq8}
\end{align}

Indeed, upon setting $b=aq/y, c=x, z=q$ and letting $a\rightarrow 0$ in this formula I obtain
\begin{align}
&(q,x)_\infty \, _2\phi_1\left(\begin{matrix*}[l] 0,0 \\ x \end{matrix*};q,q\right)
(q,y)_\infty \, _2\phi_1\left(\begin{matrix*}[l] 0,0 \\ y \end{matrix*};q,q\right)\notag \\
&=(q,x,y)_\infty \, _4\phi_3\left(\begin{matrix*}[l] \sqrt{xy/q},-\sqrt{xy/q},\sqrt{xy},-\sqrt{xy} \\
x,y,xy/q \end{matrix*};q,q\right)\label{eq9}
\end{align}
Recalling \eqref{eq6} and noting that
$$ (x)_{2n} = (-\sqrt{x},\sqrt{x},-\sqrt{xy},\sqrt{xy})_n,$$
I infer that \eqref{eq7} and \eqref{eq9} are equivalent.

In \cite{Andrews-Warnaar}, the authors utilized \eqref{eq7} with $y=-q$ to deduce that
\begin{equation}\label{eq10}
\sum_{n\geq 0}(-1)^n q^{n(n-1)/2} x^n =(x)_\infty (q;q^2)_\infty \sum_{n\geq 0}\dfrac{(-x)_{2n}q^n}{(q^2,x^2;q^2)_n}.
\end{equation}
Remarkably, if one sets $y=x/q$ in \eqref{eq3} and simplifies, one ends up with another identity for
partial theta series
\begin{equation}\label{eq11}
\sum_{n\geq 0}(-1)^n q^{n(n-1)/2} x^n =(q,x,x)_\infty \sum_{n\geq 0}\dfrac{(x^2/q^2)_{2n}q^n}{(q,x,x/q,x^2/q)_n}.
\end{equation}
Note that \eqref{eq10} and \eqref{eq11} imply a very appealing quadratic transformation formula
\begin{equation}\label{eq12}
\, _2\phi_1\left(\begin{matrix*}[l] -x, -xq \\ x^2 \end{matrix*};q^2,q\right)
=(x)_\infty(q^2;q^2)_\infty  \, _3\phi_2\left(\begin{matrix*}[l] x/\sqrt{q}, -x/\sqrt{q}, -x/q \\
x,x^2/q \end{matrix*};q,q\right).\end{equation}
To understand this result hypergeometrically, we use Heine's third transformation formula \cite[(III.3)]{Gasper-Rahman} to rewrite \eqref{eq12} as
$$
\, _2\phi_1\left(\begin{matrix*}[l] -x/q, -x \\ x^2 \end{matrix*};q^2,q^2\right)
=(x)_\infty (q;q^2)_\infty \, _3\phi_2\left(\begin{matrix*}[l] x/\sqrt{q}, -x/\sqrt{q},-x/q \\ x, x^2/q \end{matrix*};q,q\right).$$
The above can be recognized as a special case ($a=-x/q, b=x/\sqrt{q}, z=-q$) of Jane's quadratic formula
\cite[(Ex. 3.2.i)]{Gasper-Rahman}
$$\, _3\phi_2\left(\begin{matrix*}[l] a, b, -b \\ b^2, az \end{matrix*};q,-z\right)
=\dfrac{(z)_\infty}{(az)_\infty}
\, _2\phi_1\left(\begin{matrix*}[l] a, aq\\ qb^2 \end{matrix*};q^2,z^2\right).
$$
I now move on to prove 
\begin{theorem}\label{Th1}
Formula \eqref{eq4} holds true.
\end{theorem}

Below I will give two distinct proofs of this theorem. For my first proof I shall require nonterminating extension of the Sears-Carlitz formula for $\, _3\phi_2$ due to Gasper and Rahman \cite[(3.4.1)]{Gasper-Rahman}
\begin{align}
\notag
&\, _3\phi_2\left(\begin{matrix*}[l] a, b, c \\ aq/b, aq/c \end{matrix*};q,aqx/(bc)\right)
=\dfrac{(ax)_\infty}{(x)_\infty}\,  
_5\phi_4\left(\begin{matrix*}[l]  \sqrt{a}, -\sqrt{a}, \sqrt{aq}, -\sqrt{aq} \\
 aq/b, aq/c, ax, q/x \end{matrix*};q,q\right) \\
&+\dfrac{(a,aq/(bc),aqx/b, aqx/c)_\infty}{(aq/b, aq/c, aqx/(bc), 1/x)_\infty}\,
_5\phi_4\left(\begin{matrix*}[l] x\sqrt{a}, -x\sqrt{a}, x\sqrt{aq}, -x\sqrt{aq}, aqx/(bc) \\ aqx/b, aqx/c, xq, ax^2 \end{matrix*};q,q\right).\label{eq13}
\end{align}
Upon setting $b=1, a=y/x$ in this formula I see that
\begin{align}\notag
1 &= \dfrac{(y)_\infty}{(x)_\infty}\,
_4\phi_3\left(\begin{matrix*}[l]   \sqrt{y/x},-\sqrt{y/x},\sqrt{qy/x},-\sqrt{qy/x} \\ 
y, q/x, qy/x \end{matrix*};q,q\right) \\
&+\left(1-\dfrac{y}{x}\right) \dfrac{(qy)_\infty}{(1/x)_\infty} \, 
_4\phi_3\left(\begin{matrix*}[l] 
\sqrt{xy}, -\sqrt{xy}, \sqrt{xyq}, -\sqrt{xyq} \\ 
qx, qy, xy \end{matrix*};q,q\right)
\label{eq14}
\end{align}
Next, I rewrite \eqref{eq14} as
\begin{align}\notag
(q,x,q/x)_\infty &= (q,y,q/x)_\infty \,
_4\phi_3\left(\begin{matrix*}[l]   \sqrt{y/x},-\sqrt{y/x},\sqrt{qy/x},-\sqrt{qy/x} \\
y, q/x, qy/x \end{matrix*};q,q\right) \\
&+(y-x)(q,qx,qy)_\infty  \,
_4\phi_3\left(\begin{matrix*}[l]
\sqrt{xy}, -\sqrt{xy}, \sqrt{xyq}, -\sqrt{xyq} \\
qx, qy, xy \end{matrix*};q,q\right).
\label{eq15}
\end{align}
I now use \eqref{eq1} on the left of \eqref{eq15} and use \eqref{eq3} with $x$ replaced by $q/x$ on the first term on the right. This way I deduce, after simplification, that
\begin{align}\notag
&\sum_{n\geq 1} (-1)^n q^{n(n-1)/2} (x^n-y^n) \\ \notag
&= (y-x)(q,qx,qy)_\infty  \,
_4\phi_3\left(\begin{matrix*}[l]
\sqrt{xy}, -\sqrt{xy}, \sqrt{xyq}, -\sqrt{xyq} \\
qx, qy, xy \end{matrix*};q,q\right) \\
&=(y-x)(q,qx,qy)_\infty \sum_{n\geq 0} \dfrac{(xy)_{2n} q^n}{(q,qx,qy,xy)_n},
\label{eq16}
\end{align}
as desired.

To illuminate further the relation between \eqref{eq3} and \eqref{eq4} I now show how to deduce
\eqref{eq3} from \eqref{eq4} in a completely elementary, if not entirely trivial, fashion. 
To achieve this I begin with an easily verifiable identity
$$\dfrac{x^{n+1}-y^{n+1}}{x-y} = (x^n+y^n) + xy \dfrac{x^{n-1}-y^{n-1}}{x-y}.$$
This way I am led to 
\begin{equation} \label{eq17} 1+\sum_{n\geq 1} (-1)^n q^{n(n-1)/2}(x^n+y^n)=-L(x_1,y_1)
+qxy L(x_2,y_2),\end{equation}
where
$$
L(x,y) =\sum_{n\geq 1} (-1)^n q^{n(n-1)/2} \dfrac{x^n-y^n}{x-y}$$
and $x_1=x/q, y_1=y/q, x_2=qx, y_2=qy.$
Next, I employ
$$\dfrac{(x_1y_1)_{2n}}{(x_1y_1)_n} -\dfrac{(xy/q)_{2n}}{(xy)_n}
=\dfrac{xy}{q}(xyq^n)_{n-2}(1-q^n)(1-q^{n-1})$$
to infer that
\begin{align}\notag
(q,x,y)_\infty \sum_{n\geq 0} \dfrac{(xy/q)_{2n} q^n}{(q,x,y,xy)_n}
&= (q,qx_1,qy_1)_\infty \sum_{n\geq 0} \dfrac{(x_1y_1)_{2n}q^n}{(q,qx_1,qy_1,x_1y_1)_n}\\
&- xyq(q,qx_2,qy_2)_\infty \sum_{n\geq 0} \dfrac{(x_2y_2)_{2n}q^n}{(q,qx_2,qy_2,x_2y_2)_n}.
\label{eq18}
\end{align}
Combining \eqref{eq18}, \eqref{eq4} and \eqref{eq17} I arrive at \eqref{eq3}.

It is a bit more of a challenge to deduce \eqref{eq4} from \eqref{eq3}. Combining
\eqref{eq18}, \eqref{eq3} and \eqref{eq17}, one gets the following $q$-difference equation
\begin{equation}\label{eq19}
F(x,y) = q^3xy F(xq^2,yq^2),
\end{equation}
where
$$F(x,y) := L(x,y)+(q,qx,qy)\sum_{n\geq 0} \dfrac{(xy)_{2n}q^n}{(q,qx,qy,xy)_n}.$$
Iterating \eqref{eq19}, we see tht 
$$F(x,y) = 0.$$
And so, \eqref{eq4} holds.

Thus, it is possible to bypass the transformation \eqref{eq13} in establishing the equivalences of 
\eqref{eq3} and \eqref{eq4}. However, I would like to emphasize that this transformation was indispensable
in discovering \eqref{eq4} in the first place.

\section{Concluding remarks}

There are various corollaries that follow from Theorem \ref{Th1}. For example, if I set $y=-x$ and then
replace $x$ by $\sqrt{x/q}$ in \eqref{eq4}, I obtain that
\begin{equation}
\label{eq20}
\sum_{n\geq 0}q^{2n^2} x^n = (q)_\infty (xq;q^2)_\infty
\sum_{n\geq 0}\dfrac{(-x/q)_{2n} q^n}{(q,-x/q)_n(xq;q^2)_n}.
\end{equation}
Recalling Gauss' formula 
$$\sum_{n\geq 0} q^{n(n+1)/2} =\dfrac{(q^2;q^2)_\infty^2}{(q)_\infty},$$
we can deduce from \eqref{eq20} with $x=q^2$ that
$$(1-q)\dfrac{(q^8;q^8)_\infty^2 (q^2;q^2)_\infty}{(q^4;q^4)_\infty (q)_\infty^2}
=\, _2\phi_1\left(\begin{matrix*}[l] -q^2,-q \\ q^3 \end{matrix*} ; q^2,q\right).$$
Remarkably, this is a special case of the $q$-Kummer identity \cite[(11.9)]{Gasper-Rahman}.

Or, I can set $x=1$ and then let $y\rightarrow 1$ in \eqref{eq4} to obtain that
$$\sum_{n\geq 1} (-1)^n n q^{n(n-1)/2} = -(q)_\infty^3\left\{1+\sum_{n\geq 1} \dfrac{(q)_{2n-1} q^n}{(q)_n^3 (q)_{n-1}}\right\}.$$
The above can be streamlined further with the aid of \eqref{eq2}. This way I derive
$$\dfrac{1}{(q)_\infty^3} \sum_{n\geq 1} (-1)^{n+1} n q^{n(n+1)/2} = \sum_{n\geq 1} \left[\begin{matrix} 2n-1 \\
n \end{matrix}\right]_q\dfrac{q^n}{(q)_n^2},$$
where for $n,m\in\mathbf N$
$$\left[\begin{matrix} m+n \\
n \end{matrix}\right]_q :=\dfrac{(q)_{m+n}}{(q)_n (q)_m}.$$
The last identity screams for a combinatorial explanation. Perhaps, it can be found alone the lines of \cite{Berkovich-Garvan}. 

Finally, I replace $q$ by $q^4$, $x$ by $-xq^2$ in \eqref{eq10} and compare the result with \eqref{eq20}. This yields the following octonic transformation formula
\begin{align*}
&\, _2\phi_1\left(\begin{matrix*}[l] xq^2, xq^6 \\ x^2q^4 \end{matrix*};q^8,q^4\right)\\
&=\dfrac{(q)_\infty (xq;q^2)_\infty}{(q^4;q^8)_\infty (-xq^2;q^4)_\infty}\,
_4\phi_3\left(\begin{matrix*}[l] \sqrt{-x}, -\sqrt{-x},\sqrt{-x/q},-\sqrt{-x/q} \\ -x/q, \sqrt{xq}, -\sqrt{xq}\end{matrix*};q,q\right).
\end{align*}

\medskip
\bigskip

\noindent 
\textbf{Acknowledgement.} My very special thanks to Heng Huat Chan for encouraging me to submit this old work for publication and for his kindness and time taken to retype this manuscript.

\end{document}